\documentclass[12pt]{article}
\usepackage{amsmath}
\usepackage{amsfonts}
\usepackage{amscd}
\usepackage[english]{babel}
\usepackage[cp866]{inputenc}

\newtheorem{theorem}{Theorem}
\newtheorem{lemma}{Lemma}
\newtheorem{corollary}{Corollary}
\newtheorem{proposition}{Proposition}

\newtheorem{definition}{Definition}

\def\g {\mathfrak{g}}
\def\h {\mathfrak{h}}
\def\p {\mathfrak{p}}
\title{Extreme values of sectional curvature on the homogeneous
complex manifolds $U(n+1)/U(n)\times U(p+1)/U(p)$}
\author{Daurtseva N. A.,\ Kemerovo State
University}

\begin{document}
\maketitle
\begin{abstract}
Invariant complex structures on the homogeneous manifold
$U(n+1)/U(n)\times U(p+1)/U(p)$ are reseached. Extreme values of
sectional curvature of Hermitian metrics on this manifold are
found.
\end{abstract}
\centerline{\bf{1. Hermitean structures on the $S^{2n+1}\times
S^{2p+1}$}}

Recall some known construction of complex structures on
$S^{2n+1}\times S^{2p+1}\quad$ \cite{2}. It is known that
$S^{2n+1}\times S^{2p+1}$ is a principal $S^1\times S^1$ bundle
over $\mathbb{CP}^n\times\mathbb{CP}^p$. The space
$\mathbb{CP}^n\times\mathbb{CP}^p$ and fiber $S^1\times S^1$ are
complex manifolds. If we fix complex structures on the base and
fiber, then we can choose holomorphic transition functions to get
complex structure on $S^{2n+1}\times S^{2p+1}$. All those
structures form two parametric family $I(a,c)$ ($c>0$).

The product $S^{2n+1}\times S^{2p+1}$ is homogeneous space
$U(n+1)/U(n)\times U(p+1)/U(p)$. Let $\g_1$ and $\h_1$ ($\g_2$ and
$\h_2$) are Lee algebras $\mathfrak{u}(n+1)$ and $\mathfrak{u}(n)$
( $\mathfrak{u}(p+1)$ and $\mathfrak{u}(p)$) of Lee groups
$U(n+1)$ and $U(n)$ ($U(p+1)$ and $U(p)$). Let
$\h=\h_1\times\h_2$. The complex structures $I(a,c)$ are $ad\h$ -
invariant \cite{2}.

Let $E^1_{\nu\mu}$ is matrix with 1 on the $(\nu,\mu)$ - place and
other zero elements. Define:
$$
Z^1_{\nu\mu}=E^1_{\nu\mu}-E^1_{\mu\nu},\
T^1_{\nu\mu}=E^1_{\nu\mu}+E^1_{\mu\nu},\ 0\leq\mu <\nu\leq n,
$$
Take decomposition $\g_1=\h_1\oplus \p_1$, where $\mathfrak{p}_1$
has basis $X^1=\frac12iT^1_{00}$, $Y^1_{2\nu-1}=Z^1_{\nu 0}$,
$Y^1_{2\nu}=iT^1_{\nu 0}$. Remark that $X^1$ is tangent to fiber
$S^1$ in the fibre bundle
$S^{2n+1}\longrightarrow^{S^1}\mathbb{CP}^n$. By analogy, we have
$\g_2=\h_2\oplus\p_2$ and $\g_1\times\g_2=\h\oplus\p$, where
$\p=\p_1\times\p_2$. So $S^{2n+1}\times S^{2p+1}$ viewed as
homogeneous space $U(n+1)/U(n)\times U(p+1)/U(p)$ has basis
$X^1,Y^1_{2\nu-1},Y^1_{2\nu},X^2,Y^2_{2\mu-1},Y^2_{2\mu}$,
$1\leq\nu\leq n$, $1\leq\mu\leq p$. On these basis vectors
$$
I(a,c)X^1=\frac{a}{c}X^1+\frac1cX^2,\quad
I(a,c)X^2=-\frac{a^2+c^2}{c}X^1-\frac{a}{c}X^2
$$
$$
I(a,c)Y^1_{2\nu-1}=Y^1_{2\nu},\quad I(a,c)Y^2_{2\mu-1}=Y^2_{2\mu},
$$
parameters $a$ and $c$ are real and $c>0$.
\begin{definition}
Almost complex structure $J$ on the manifold $M$ is called
positive associated with 2-form $\omega$ if:\\
1) $\omega(JX,JY)=\omega(X,Y),$, for all $X,Y\in TM$\\
2) $\omega(X,JX)>0$, for all nonzero $X\in TM$
\end{definition}
Fix on the $S^{2n+1}\times S^{2p+1}$ nondegenerate invariant
2-form $\omega$  by:
$$
\omega=X^1\wedge X^2+\sum_{\nu = 1}^nY^1_{2\nu-1}\wedge
Y^1_{2\nu}+\sum_{\nu = 1}^pY^2_{2\nu-1}\wedge Y^2_{2\nu}
$$
\begin{lemma}
All complex structures $I(a,c)$ are positive associated with
$\omega$.
\end{lemma}
{\bf Proof.} For $I(a,c)$ properties 1) and 2) of definition 1 are
obvious
\begin{corollary}
Each complex structure $I(a,c)$ defines unique $\omega$ -
associated metric by formula
$$
g(a,c)(X,Y)=\omega(X,I(a,c)Y)
$$
\end{corollary}
These associated metrics in the above basis are:
$$
g(a,c)=\left(\begin{array}{cc}
g_{11} & g_{12}\\
g_{21} & g_{22}
\end{array}\right),
$$
where
$$
g_{11}=\left(
\begin{array}{cccc}
1/c & 0^{\top}_n\\
0_n & E_n
\end{array} \right),\ g_{22}= \left(
\begin{array}{cccc}
(a^2+c^2)/c & 0^{\top}_p\\
0_p & E_p
\end{array} \right),
$$
$$
g^{\top}_{21}=g_{12}=\left(\begin{array}{ccc} -a/c & 0^{\top}_p\\
0_n & 0\end{array}\right)
$$
$0_n$ is zero column-vector, $E_n$ is unit $n\times n$-matrix.

Each metric of this family is $I(a,c)$-Hermitean, so we obtain
two-parametric family of Hermitean manifolds $(S^{2n+1}\times
S^{2p+1},$ $g(a,c), I(a,c), \omega)$. Invariant metric induces
scalar product on $\p$. We will denote this product as
$\langle,\rangle_{a,c}$.
\begin{proposition}
Invariant Riemannian connection for $g(a,c)$ on $S^{2n+1}\times
S^{2p+1}$ is given by formula $D_XY=\frac12[X,Y]_{\p}+U(X,Y)$,
where $U$ is symmetric bilinear mapping
$U:\p\times\p\longrightarrow\p$:
$$
U(X^1,Y^1_{2\nu-1})=\frac{2-c}{2c}\ Y^1_{2\nu},\
U(X^1,Y^1_{2\nu})=-\frac{2-c}{2c}\ Y^1_{2\nu-1},
$$
$$
U(X^1,Y^2_{2\nu-1})=-\frac{a}{c}\ Y^2_{2\nu},\
U(X^1,Y^2_{2\nu})=\frac{a}{c}\ Y^2_{2\nu-1},
$$
$$
U(X^2,Y^1_{2\nu-1})=-\frac{a}{c}\ Y^1_{2\nu},\
U(X^2,Y^1_{2\nu})=\frac{a}{c}\ Y^1_{2\nu-1},
$$
$$
U(X^2,Y^2_{2\nu-1})=\left(\frac{a^2+c^2}{c}-\frac12\right)\
Y^2_{2\nu},\
U(X^2,Y^2_{2\nu})=-\left(\frac{a^2+c^2}{c}-\frac12\right)\
Y^2_{2\nu-1},
$$
For other basis vectors $X$ and $Y$ the $U(X,Y)$ is equal to 0.
\end{proposition}
{\bf Proof.} One can find $U$ by formula:
$2g(U(X,Y),Z)=g([Z,X]_{\p},Y)+g(X,[Z,Y]_{\p})$

\begin{proposition}
Two-parametric family of metrics $g(a,c)$ has following
characteristics:\\
1) Ricci curvature in the above basis is:
$$Ric(a,c)=\left(\begin{array}{cc}
r_{11} & r_{12}\\
r_{21} & r_{22}
\end{array}\right), r_{11}=\left(\begin{array}{cccc} 2\frac{n+pa^2}{c^2} & 0^{\top}_n\\
0_n & 2(1+n-\frac{1}{c})E_n
\end{array}\right),
$$
$$
r_{22}=\left(\begin{array}{cccc} 2\frac{na^2+p(a^2+c^2)^2}{c^2} & 0^{\top}_p\\
0_p & 2(1+p-\frac{a^2+c^2}{c})E_p
\end{array}\right),
$$
$$
r^{\top}_{21}=r_{12}=\left(\begin{array}{ccc} -2\frac{a}{c^2}(n+p(a^2+c^2)) & 0^{\top}_p\\
0_n & 0\end{array}\right).
$$
Proper values of Ricci curvature $\tilde{r}_i$ are
$\tilde{r}_{1,2}=\frac{x+y\pm\sqrt{(x-y^2+4z^2)}}{2}$, where
$x=2\frac{n+pa^2}{c^2}$, $y=2\frac{na^2+p(a^2+c^2)^2}{c^2}$,
$z=-2\frac{a}{c^2}(n+p(a^2+c^2))$;
$\tilde{r}_3=\tilde{r}_4=\dots=\tilde{r}_{2n+2}=2(1+n-\frac1c)$,
$\tilde{r}_{2n+3}=\tilde{r}_{2n+4}=\dots=\tilde{r}_{2n+2p+2}=2(1+p-\frac{a^2+c^2}{c})$.\\
2). Scalar curvature is given by following formula:
$$
s=4n\left(1+n-\frac{1}{2c}\right)+4p\left(1+p-\frac{a^2+c^2}{2c}\right)
$$
\end{proposition}
{\bf Proof.} Compute Ricci curvature by formula \cite{1}
$$Ric(a,c)(X,X)=-\frac12\sum_i|[X,v_i]_{\p}|^2-\frac12\sum_i\langle[X,[X,v_i]_{\p}]_{\p},v_i\rangle_{a,c}$$
$$-\sum_i\langle[X,[X,v_i]_{\h}]_{\p},v_i\rangle_{a,c}+\frac14\sum_{i,j}\langle[v_i,v_j]_{\p},X\rangle_{a,c}^2-\langle[Z,X]_{\p},X\rangle_{a,c}$$
where $Z=\sum_iU(v_i,v_i)$ and $v_i$ is orthonormal basis of
$(\p,\langle,\rangle_{a,c})$. Scalar curvature is computed as a
trace of Ricci curvature: $s=Ric_{ij}g^{ij}$, where $g^{ij}$ are
components of $g(a,c)^{-1}$ ($i,j=1,\dots,2n+2p+2$).
\newpage
\centerline{\bf{2. Extreme values of sectional curvature of
metrics $g(a,c)$.}}

Use bivectors to compute sectional curvature. Let $X_0,\dots,X_m$
is orthonormal basis of tangent bundle of any manifold. Let
$B_{\nu\mu}=X_{\nu}\wedge X_{\mu}$ is basis of space of bivectors.

Let define:
$$
R_{(\alpha\nu)(\rho\mu)}=\langle
R(X_{\alpha},X_{\nu})X_{\mu},X_{\rho}\rangle_{a,c}
$$
$R$ will denote $\frac{m(m+1)}{2}\times\frac{m(m+1)}{2}$ - matrix
($\frac{m(m+1)}{2}$ is dimension of space of bivectors) with
element $R_{(\alpha\nu)(\rho\mu)}$ on the place
$((\alpha\nu)(\rho\mu))$, where $\alpha<\nu$, and $\rho<\mu$.

As curvature is polylinear we have: sectional curvature in
direction of unit bivector $B$ is equal to
$$
K(B)=\langle R\cdot B^T,B\rangle_{a,c},
$$
where $B^T$ is column-vector, which consists of bivector $B$
coordinates, $R\cdot B^T$ is matrix-product,
$\langle,\rangle_{a,c}$ is scalar product in the space of
bivectors, corresponding to metric $\langle,\rangle_{a,c}$ on
vectors.

Bivector $B$ is called decomposable, if there exist vectors $X$
and $Y$ such as $B=X\wedge Y$. Note, that only decomposable
bivectors correspond to 2-dimensional plane in the tangent bundle.

Compute matrix of curvature operator $R$ of space $(S^{2n+1}\times
S^{2p+1}, g(a,c))$. Orthonormar basis of space
$(\p,\langle,\rangle_{a,c})$ is
$$
Z_0=\sqrt{c}X^1,\ Z_{\nu}=Y^1_{\nu},\ \nu=1,\dots,2n
$$
$$
Z_{2n+1}=\frac{a}{\sqrt{c}}X^1+\frac{1}{\sqrt{c}}X^2,\
Z_{2n+1+\mu}=Y^2_{\mu},\ \mu=1,\dots,2p
$$
\begin{proposition}
$$
R_{(0\nu)(\rho\mu)}=\left\{
\begin{array}{l}
\frac1c,\ \rho=0,\ \nu=\mu=1,\dots,2n,\\
\frac{a^2}{c},\ \rho=0,\ \nu=\mu=2n+2,\dots,2n+2p+1,\\
-a,\ \rho=2n+1,\ \nu=\mu=2n+2,\dots,2n+2p+1,\\
0,\mbox{ for other }\nu,\rho,\mu.
\end{array}
\right.
$$
$$
R_{(2n+1\nu)(\rho\mu)}=\left\{
\begin{array}{l}
-a,\ \rho=0,\ \nu=\mu=2n+2,\dots,2n+2p+1,\\
c,\ \rho=2n+1,\ \nu=\mu=2n+2,\dots,2n+2p+1,\\
0,\mbox{ for other }\nu,\rho,\mu.
\end{array}
\right.
$$
\end{proposition}
{\bf Proof.} Compute $R_{(0\nu)(\rho\mu)}$. By definition
$$
R_{(0\nu)(\rho\mu)}=-\sqrt{c}\langle
R(X^1,Z_{\nu})Z_{\rho},Z_{\mu}\rangle_{a,c}.
$$
$$
R(X^1,Z_{\nu})Z_{\rho}=D_{X^1}\left(\frac12[Z_{\nu},Z_{\rho}]_{\p}+U(Z_{\nu},Z_{\rho})\right)-
D_{Z_{\nu}}\left(\frac12[X^1,Z_{\rho}]_{\p}+U(X^1,Z_{\rho})\right)
$$
$$
-\left(\frac12[[X^1,Z_{\nu}]_{\p},Z_{\rho}]_{\p}+U([X^1,Z_{\nu}]_{\p},Z_{\rho})\right)
$$
1) Let $\nu=1,\dots,2n$, study the case of $\nu=2k$, nonzero
values of $R$ are possible in two cases only,when $\rho=2k$ or
$\rho=2n+1$.
\\
-- $\rho=2k$:
$$
R(X^1,Y^1_{2k})Y^1_{2k}=-D_{Y^1_{2k}}\left(1-\frac1c\right)Y^1_{2k-1}+X^1=\frac1cX_1
$$
-- $\rho=2n+1$:
$$
R(X^1,Y^1_{2k})Z_{2n+1}=-D_{X^1}\left(-\frac{a}{2\sqrt{c}}+\frac{a(c-2)}{2c\sqrt{c}}+\frac{a}{c\sqrt{c}}\right)Y^1_{2k-1}
$$
$$
-\left(\frac{a}{2\sqrt{c}}+\frac{a(2-c)}{2c\sqrt{c}}-\frac{a}{c\sqrt{c}}\right)Y^1_{2k}=0
$$
By analogy, for $\nu=2k-1$,
$$
R(X^1,Y^1_{2k-1})Z_{\rho}=\left\{\begin{array}{ll} \frac{1}{c}X^1,
& \rho=2k-1,\\
0, & \rho\neq 2k-1.
\end{array}\right.
$$
So, for $\nu=1,\dots,2n$, $\rho\neq 1$:
$$
R_{(0\nu)(\rho\mu)}=\left\{\begin{array}{ll} \frac1c, & \rho=0,\
\nu=\mu,\\
0, & \mbox{in the otner cases}.\end{array}\right.
$$
2) Let $\nu=2n+1$, we have:
$$
R(X^1,Z_{2n+1})Z_{\rho}=\frac{a}{\sqrt{c}}R(X^1,X^1)Z_{\rho}+\frac{1}{\sqrt{c}}R(X^1,X^2)Z_{\rho}=0
$$
3) Let $\nu=2n+1+i$, where $i=1,\dots,2p$. Study the case
$i=2k-1$, we have:
$$
R(X^1,Y^2_{2k-1})Z_{\rho}=D_{X^1}(\frac12[Y^2_{2k-1},Z_{\rho}]_{\p}+U(Y^2_{2k-1},Z_{\rho}))
-D_{Y^2_{2k-1}}(\frac12[X^1,Z_{\rho}]_{\p}+U(X^1,Z_{\rho}))
$$
Nonzero values of $R$ are possible in cases $\rho=2n+1$ or
$\rho=2n+2k$:\\
-- $\rho=2n+1$:
$$
R(X^1,Y^2_{2k-1})Z_{2n+1}=D_{X^1}(\frac12[Y^2_{2k-1},\frac{1}{\sqrt{c}}X^2]_{\p}+
U(Y^2_{2k-1},\frac{a}{\sqrt{c}}X^1+\frac{1}{\sqrt{c}}X^2))=\frac{a}{\sqrt{c}}Y^2_{2k-1}
$$
-- $\rho=2n+2k$:
$$
R(X^1,Y^2_{2k-1})Z_{2n+2k}=-D_{Y^2_{2k-1}}U(X^1,Y^2_{2k-1})=-\frac{a}{c}X^2
$$
By analogy, for $i=2k$ we have:
$$
R(X^1,Y^2_{2k})Z_{2n+1}=\frac{a}{\sqrt{c}}Y^2_{2k},\qquad
R(X^1,Y^2_{2k})Z_{2n+2k+1}=-\frac{a}{c}Y^2_{2k}
$$
So, for $\nu=2n+2,\dots,2n+2p+1$
$$
R_{(0\nu)(\rho\mu)}=\left\{\begin{array}{ll} \frac{a^2}{c}, &
\rho=0,\ \nu=\mu,\\
-a, &
\rho=2n+1,\ \mu=\nu,\\
0, & \rho\neq 2n+1,\ \rho\neq 0.
\end{array}\right.
$$
Compute $R_{(2n+1\nu)(\rho\mu)}$. By definition:
$$
R_{(2n+1\nu)(\rho\mu)}=-\langle
R(Z_{2n+1}Z_{\nu})Z_{\rho},Z_{\mu}\rangle_{a,c}=-\frac{a}{\sqrt{c}}\langle
R(X^1,Z_{\nu})Z_{\rho},Z_{\mu}\rangle_{a,c}
-\frac1{\sqrt{c}}\langle
R(X^2,Z_{\nu})Z_{\rho},Z_{\mu}\rangle_{a,c}
$$
Let $\nu=2n+1+i$, as:
$$
R_{(2n+1\nu)(\rho\mu)}=R_{(\rho\mu)(2n+1\nu)}
$$
it is enough to research case $\rho>0$. \\
1) $\rho=1,\dots,2n$
$$
R(X^1,Y^2_i)Y^1_{\rho}=-D_{Y^2_i}\left(\frac12[X^1,Y^1_{\rho}]_{\p}+U(X^1,Y^1_{\rho})\right)=0
$$
$$
R(X^2,Y^2_i)Y^1_{\rho}=-D_{Y^2_i}U(X^2,Y^1_{\rho})=0
$$
2) $\rho=2n+1$\\
Values of $R(X^1,Y^2_i)Z_{2n+1}$ are obtained above.
$$
R(X^2,Y^2_i)Z_{2n+1}=\frac{a}{\sqrt{c}}R(X^2,Y^2_i)X^1+\frac{1}{\sqrt{c}}R(X^2,Y^2_i)X^2
$$
$$
R(X^2,Y^2_i)X^1=D_{X^2}U(Y^2_i,X^1)-U([X^2,Y^2_i]_{\p},X^1)=
$$
$$
=\left\{\begin{array}{ll} i=2k-1: &
-\frac{a}{c}D_{X^2}Y^2_{2k}+U(Y^2_{2k},X^1)=\frac{a(a^2+c^2)}{c^2}Y^2_{2k-1},\\
i=2k: &
\frac{a}{c}D_{X^2}Y^2_{2k-1}-U(Y^2_{2k-1},X^1)=\frac{a(a^2+c^2)}{c^2}Y^2_{2k}.
\end{array}\right.
$$
$$
R(X^2,Y^2_i)X^2=D_{X^2}\left(\frac12[Y^2_i,X^2]_{\p}+U(Y^2_i,X^2)\right)
-\frac12[[X^2,Y^2_i]_{\p},X^2]_{\p}-U([X^2,Y^2_i]_{\p},X^2)
$$
$$
=\left\{\begin{array}{ll} i=2k-1: &
\frac{a^2+c^2}{c}D_{X^2}Y^2_{2k}-\frac{a^2+c^2}{c}Y^2_{2k-1}=-\frac{(a^2+c^2)^2}{c^2}Y^2_{2k-1},\\
i=2k: &
-\frac{a^2+c^2}{c}D_{X^2}Y^2_{2k-1}-\frac{a^2+c^2}{c}Y^2_{2k}=-\frac{(a^2+c^2)^2}{c^2}Y^2_{2k}.
\end{array}\right.
$$
We obtain:
$$
R(X^2,Y^2_i)Z_{2n+1}=-\frac{a^2+c^2}{\sqrt{c}}Y^2_i
$$
We have:
$$
R_{(2n+1\nu)(2n+1\mu)}=-\frac{a}{\sqrt{c}}\left(\frac{a}{\sqrt{c}}Y^2_i,Z_{\mu}\right)
-\frac1{\sqrt{c}}\left(-\frac{a^2+c^2}{\sqrt{c}}Y^2_i,Z_{\mu}\right)
$$
$$
=\left\{\begin{array}{ll} c, &  \mu=2n+1+i,  \\
0, & \mu\neq 2n+1+i.
\end{array}\right.
$$
3) Let $\rho=2n+1+j$
$$
R(X^1,Y^2_i)Y^2_j=-D_{Y^2_i}U(X^1,Y^2_j)=\left\{\begin{array}{ll}
-\frac{a}{c}X^2, & i=j,\\
0, & i\neq j.
\end{array}\right.
$$
$$
R(X^2,Y^2_i)Y^2_j=-D_{Y^2_i}\left(\frac12[X^2,Y^2_j]_{\p}+U(X^2,Y^2_j)\right)-\frac12[[X^2,Y^2_i]_{\p},Y^2_j]_{\p}
-U([X^2,Y^2_i]_{\p},Y^2_j)
$$
$$
=\left\{\begin{array}{ll} \frac{a^2+c^2}{c}X^2, & i=j,\\
0, & i\neq j.\end{array}\right.
$$
So, for $\rho=2n+2,\dots,2n+2p+1$ $R_{(2n+1\nu)(\rho\mu)}=0$, if
$\mu>\rho$.

Let $B$ is some decomposable bivector, having coordinates
$(b_{01},\dots,b_{2n+2p2n+2p+1})$ in the basis $Z_{\nu}\wedge
Z_{\mu}$. Decompose $B$ in the following sum:
$$
B=B'+B''
$$
where\\
$B'=(b_{01},\dots,b_{02n},0,b_{02n+2},\dots,b_{02n+2p+1},0,\dots,0,b_{2n+12n+2},\dots,b_{2n+12n+2p+1},0,\dots,0)$,
$B''=B-B'$. By proposition 3:
$$
K(B)=\langle R\cdot B^T,B\rangle_{a,c}=\langle R\cdot
B'^T,B'\rangle_{a,c}+\langle R\cdot B''^T,B''\rangle_{a,c}
$$
If decompose $B'$ in the sum:
$$
B'=B'_1+B'_2+B'_3,
$$
where\\
$B'_1=(b_{01},\dots,b_{02n},0,\dots,0)$,\\
$B'_2=(0,\dots,0,b_{02n+2},\dots,b_{02n+2p+1},0,\dots,0)$,\\
$B'_3=(0,\dots,0,b_{2n+12n+2},\dots,b_{2n+12n+2p+1},0,\dots,0)$,
then:
$$
K(B)=\frac1c\|B'_1\|^2+\frac{a^2}{c}\|B'_2\|^2+c\|B'_3\|^2-2a\sum_{i=1}^pb_{02n+1+i}b_{2n+12n+1+i}+\langle
R\cdot B''^T,B''\rangle_{a,c}
$$
\begin{proposition}
$R_{(2l-12l)(2n+2m2n+2m+1)}=\frac{2a}{c}$\\
$R_{(\alpha\nu)(\rho\mu)}=0$ for other $\alpha,\nu=1,\dots,2n$,
$\rho,\mu=2n+2,\dots,2n+2p+1$
\end{proposition}
{\bf Proof.} Let $\alpha,\nu=1,\dots,2n$, $\rho=2n+1+i$,
$\mu=2n+1+j$, $i,j=1,\dots,2p$, then:
$$
R_{(\alpha\nu)(\rho\mu)}=\langle
R(Y^1_{\alpha},Y^1_{\nu})Y^2_j,Y^2_i\rangle_{a,c}
$$
$$
R(Y^1_{\alpha},Y^1_{\nu})Y^2_j=-\frac12[[Y^1_{\alpha},Y^1_{\nu}]_{\p},Y^2_j]_{\p}
-U([Y^1_{\alpha},Y^1_{\nu}]_{\p},Y^2_j)=
$$
$$
=\left\{\begin{array}{ll} 2U(X^1,Y^2_j) & \alpha=2l-1,\ \nu=2l,\\
0, & \mbox{for other }\alpha,\ \nu.
\end{array}\right.
$$
As
$$
2U(X^1,Y^2_j)=\left\{\begin{array}{ll}
\frac{2a}{c}Y^2_{2m-1}, & j=2m,\\
-\frac{2a}{c}Y^2_{2m}, & j=2m-1.
\end{array}\right.,
$$
then
$$
R_{(2l-12l)(2n+2m2n+2m+1)}=\frac{2a}{c}
$$

\begin{proposition}
$$
R_{(2l-12n+2m)(2l2n+2m+1)}=\frac{a}{c}
$$
$$
R_{(2l2n+2m)(2l-12n+2m+1)}=-\frac{a}{c}
$$
$$
R_{(2l-12n+2m+1)(2l2n+2m)}=-\frac{a}{c}
$$
$$
R_{(2l2n+2m+1)(2l-12n+2m)}=\frac{a}{c}
$$
$$
R_{(\alpha\nu)(\rho\mu)}=0
$$
for other values of $\alpha,\rho=1,\dots,2n$,
$\nu,\mu=2n+2,\dots,2n+2p+1$
\end{proposition}
{\bf Proof.} Let $\alpha,\rho=1,\dots,2n$, $\nu=2n+1+i$,
$\mu=2n+1+j$.
$$
R_{(\alpha\nu)(\rho\mu)}=\langle
R(Y^1_{\alpha},Y^2_i)Y^2_j,Y^1_{\rho}\rangle_{a,c}
$$
$$
R(Y^1_{\alpha},Y^2_i)Y^2_j=\frac12D_{Y^1_{\alpha}}[Y^2_i,Y^2_j]_{\p}=\left\{\begin{array}{ll}
-D_{Y^1_{\alpha}}X^2, & i=2m-1,\ j=2m,\\
D_{Y^1_{\alpha}}X^2, & i=2m,\ j=2m-1,\\
0, & |j-i|\neq 1\end{array}\right.
$$
$$
=\left\{\begin{array}{ll} -U(Y^1_{2l-1},X^2)=\frac{a}{c}Y^1_{2l},
&
\alpha=2l-1,\ i=2m-1,\ j=2m,\\
-U(Y^1_{2l},X^2)=-\frac{a}{c}Y^1_{2l}, &
\alpha=2l,\ i=2m-1,\ j=2m,\\
U(Y^1_{2l-1},X^2)=-\frac{a}{c}Y^1_{2l}, &
\alpha=2l-1,\ i=2m,\ j=2m-1,\\
U(Y^1_{2l},X^2)=\frac{a}{c}Y^1_{2l}, & \alpha=2l,\ i=2m,\ j=2m-1,\\
0, & |j-i|\neq 1.
\end{array}\right.
$$

Decompose $B''$ in sum of three bivectors:
$$
B''=B_1''+B_2''+B_3''
$$
where $B_1''$ is bivector with nonzero coordinates
$b_{ij}$, $i=1,\dots,2n$, $j=2n+2,\dots,2n+2p+1$,\\
$B''_2$ is bivector with nonzero coordinates $b_{ij}$,
$i,j=1,\dots,2n$,\\
$B''_3$ is bivector with nonzero coordinates $b_{ij}$,
$i,j=2n+2,\dots,2n+2p+1$.\\
Then by results of propositions 4 and 5 we have:
$$
\langle R\cdot B''^T,B''\rangle_{a,c}=\langle R\cdot
B''^T_1,B''_1\rangle_{a,c}+\langle R\cdot
B''^T_2,B''_3\rangle_{a,c}+\langle R\cdot
B''^T_3,B''_2\rangle_{a,c}+\langle R\cdot
B''^T_2,B''_2\rangle_{a,c}
$$
$$
+\langle R\cdot
B''^T_3,B''_3\rangle_{a,c}=2\frac{a}{c}\sum_{l=1}^n\sum_{m=1}^p(b_{2l-12n+2m}b_{2l2n+2m+1}-b_{2l2n+2m}b_{2l-12n+2m+1})
$$
$$
+4\frac{a}{c}\sum_{l=1}^n\sum_{m=1}^pb_{2l-12l}b_{2n+2m2n+2m+1}+\langle
R\cdot B''^T_2,B''_2\rangle_{a,c}+\langle R\cdot
B''^T_3,B''_3\rangle_{a,c}
$$
As $\langle R\cdot B''^T_2,B''_2\rangle_{a,c}=K(B''_2)$ and
$\langle R\cdot B''^T_3,B''_3\rangle_{a,c}=K(B''_3)$ we can use
the results of work of Volper D.E. \cite{3} and obtain:
$$
\min(4-\frac{3}{c};1)\|B_2''\|^2\leq
K(B_2'')\leq\max(4-\frac{3}{c};1)\|B_2''\|^2
$$
$$
\min(4-\frac{3(a^2+c^2)}{c};1)\|B_3''\|^2\leq
K(B_3'')\leq\max(4-\frac{3}{c};1)\|B_3''\|^2
$$
\begin{theorem}
Sectional curvature $K(a,c)$ of metric $g(a,c)$ on $S^{2n+1}\times
S^{2p+1}$ satisfies to the following inequality:
$$
K_{\min}\leq K(a,c)\leq K_{\max}
$$
where
$$
K_{\min}=\left\{\begin{array}{ll}
\min(-|\frac{a}{c}|,\frac{5c-3-\sqrt{16a^2-18c+9c^2+9}}{2c})
& a^2+(c-\frac12)^2\leq\frac14\\
\min(-|\frac{a}{c}|,\frac{8c-3(1+a^2+c^2)-\sqrt{9(a^2+c^2-1)^2+16a^2}}{2c})
& a^2+(c-\frac12)^2>\frac14, c<1\\
\min(-|\frac{a}{c}|,\frac{5c-3(a^2+c^2)-\sqrt{16a^2+9c^2-18c(a^2+c^2)+9(a^2+c^2)^2}}{2c})
& c\geq 1\end{array}\right.
$$
$$
K_{\max}=\left\{\begin{array}{ll}
\max(\frac1c,\frac{5c-3(a^2+c^2)+\sqrt{16a^2+9c^2-18c(a^2+c^2)+9(a^2+c^2)^2}}{2c})
& a^2+(c-\frac12)^2\leq\frac14\\
\max(\frac{a^2+c^2}{c},1+2|\frac{a}{c}|,\frac1c) & a^2+(c-\frac12)^2>\frac14, c<1\\
\max(\frac{a^2+c^2}{c},\frac{5c-3+\sqrt{16a^2-18c+9c^2+9}}{2c}) &
c\geq 1
\end{array}\right.
$$
\end{theorem}
{\bf Proof.} As above results we have
$$
K'\leq K(a,c)\leq K''
$$
where
$$
K'=\min(4-\frac3c,1)\|B''_2\|^2+\min(4-3\frac{a^2+c^2}{c},1)\|B''_3\|^2+\frac{1}{c}\|B_1'\|^2+
\frac{a^2}{c}\|B'_2\|^2+c\|B'_3\|^2
$$
$$
-2a\sum_{m=1}^pb_{2n+12n+1+m}b_{02n+1+m}+4\frac{a}{c}\sum_{l=1}^n\sum_{m=1}^pb_{2l-12l}b_{2n+2m2n+2m+1}
$$
$$
+2\frac{a}{c}\sum_{l=1}^n\sum_{m=1}^p(b_{2l-12n+2m}b_{2l2n+2m+1}-b_{2l2n+2m}b_{2l-12n+2m+1})
$$
$$
K''=\max(4-\frac3c,1)\|B''_2\|^2+\max(4-3\frac{a^2+c^2}{c},1)\|B''_3\|^2+\frac{1}{c}\|B_1'\|^2+
\frac{a^2}{c}\|B'_2\|^2+c\|B'_3\|^2
$$
$$
-2a\sum_{m=1}^pb_{2n+12n+1+m}b_{02n+1+m}+4\frac{a}{c}\sum_{l=1}^n\sum_{m=1}^pb_{2l-12l}b_{2n+2m2n+2m+1}
$$
$$
+2\frac{a}{c}\sum_{l=1}^n\sum_{m=1}^p(b_{2l-12n+2m}b_{2l2n+2m+1}-b_{2l2n+2m}b_{2l-12n+2m+1})
$$
It is obvious that
$$
\min(4-\frac{3}{c},1)=\left\{\begin{array}{ll} 4-\frac{3}{c}, &
c<1,\\
1, & c\geq1.
\end{array}\right.
$$
$$
\min(4-3\frac{a^2+c^2}{c},1)=\left\{\begin{array}{ll} 1, &
a^2+(c-\frac12)^2<\frac14,\\
4-3\frac{a^2+c^2}{c}, & a^2+(c-\frac12)^2\geq\frac14.
\end{array}\right.
$$
Let, for example, $a^2+(c-\frac12)^2\leq\frac14$. Then
$\min(4-\frac{3}{c},1)=4-\frac{3}{c}$,
$\min(4-3\frac{a^2+c^2}{c},1)=1$. As $\|B\|=1$, we have to solve
problem of conditional extremum of function $K'(B)$. We obtain:
$$
K_{\min}=\min_{\|B\|=1}K'(B)=\min(\pm\frac{a}{c},0,\frac{a^2+c^2}{c},\frac1c,\frac{5c-3\pm\sqrt{16a^2-18c+9c^2+9}}{2c},4-3\frac1c,1)
$$
$$
=\min(-|\frac{a}{c}|,\frac{5c-3-\sqrt{16a^2-18c+9c^2+9}}{2c})
$$
By analogy
$$
K_{\max}=\max_{\|B\|=1}K''(B)
$$
$$
=\max(\pm\frac{a}{c},0,\frac{a^2+c^2}{c},\frac1c,\frac{5c-3(a^2+c^2)\pm\sqrt{16a^2+9c^2-18c(a^2+c^2)+9(a^2+c^2)^2}}{2c},
$$
$$
4-3\frac{a^2+c^2}{c},1)=\max(\frac1c,\frac{5c-3(a^2+c^2)+\sqrt{16a^2+9c^2-18c(a^2+c^2)+9(a^2+c^2)^2}}{2c})
$$

{\sc e-mail:} natali0112@ngs.ru
\end{document}